\documentclass{llncs} 
\usepackage{amssymb, amsmath}
\usepackage[all]{xypic}
\usepackage{graphicx}

\begin{document}

\title{The smallest multistationary mass-preserving chemical reaction network}

\author{Anne Shiu}

\institute{Dept. of Mathematics, University of California, Berkeley, CA 94720-3840, USA}

\maketitle

\begin{abstract}
Biochemical models that exhibit bistability are of interest to biologists and mathematicians alike.  Chemical reaction network theory can provide sufficient conditions for the existence of bistability, and on the other hand can rule out the possibility of multiple steady states.  Understanding small networks is important because the existence of multiple steady states in a subnetwork of a biochemical model
can sometimes be lifted to establish multistationarity in the larger network.  This paper establishes the smallest reversible, mass-preserving network that admits bistability and determines the semi-algebraic set of parameters for which more than one steady state exists.

\end{abstract}


\section{Introduction}
Bistable biochemical models are often presented as the possible underpinnings of chemical switches \cite{HowSwitch,Laurent}.  Systematic study of mass-action kinetics models--which {\em a priori} may or may not admit multiple steady states--constitutes chemical reaction network theory (CRNT), a subject pioneered by Horn, Jackson, and Feinberg \cite{FeinLectures,GenMAK}.  Certain classes of networks, such as those of deficiency zero, do not exhibit multistationarity or other strange behaviors.  A generalization of deficiency-zero systems is the class of {\em toric dynamical systems} which have a unique steady state \cite{TDS}.  See also the recent work of Craciun and Feinberg for additional conditions that rule out multistationarity \cite{ME1,ME2}.

On the other hand, there are conditions that are sufficient for establishing whether a network supports multiple steady states.  The CRNT Toolbox developed by Feinberg and improved by Ellison implements the Deficiency One and Advanced Deficiency Algorithms \cite{EllisonThesis,FeinDefOne}; this software is available online \cite{Toolbox}.  For a large class of systems, the CRNT Toolbox either provides a witness for multistationarity or concludes that it is impossible.  For systems for which the CRNT Toolbox is inconclusive, see the approach of Conradi {\em et al.} \cite{subnetwork}.
Related work includes an algebraic approach that determines the full set of parameters for which a system is multistationary; a necessary and sufficient condition for multistationarity is the existence of a non-trivial sign vector in the intersection of two subsets of Euclidean space \cite{MAPK}.

To model biological processes, one typically reverse-engineers a system of non-linear differential equations that exhibits specific dynamical behavior, such as bistability or oscillations, observed in experiments.  For example, Segel 
proposes a small immune network consisting only two cell types, which has three stable steady states, corresponding to ``normal,'' ``vaccinated,'' and ``diseased'' states \cite{Segel}.  Similarly, the Brusselator is a mass-action kinetics network with a stable limit cycle \cite{FeinOsc}.

This paper focuses on the smallest mass-action kinetics networks that admit multiple steady states.  Section 2 provides an introduction to chemical reaction network theory.  A special network called the {\em Square} is shown in Section 3 to be a smallest reversible multistationary chemical reaction network.  Sections 4 and 5 determine precisely which parameters of the Square give rise to multiple steady states.  

\section{Chemical Reaction Network Theory}

We now give an introduction to chemical reaction network theory.  Before giving precise definitions, we present an intuitive example that illustrates how a chemical reaction network gives rise to a dynamical system.  An example of a {\em chemical reaction}, as it usually appears in the literature, is the following:
\[
 \begin{xy}<15mm,0cm>:
 (1,0) ="3A+C" *+!L{3A+C} *{};
 (0,0) ="A+B" *+!R{A+B} *{};
 (0.55,.05)="k" *+!D{\kappa} *{};
   {\ar "A+B"*{};"3A+C"*{}};
    \end{xy}\]
\\
In this reaction, one unit of chemical {\em species} $A$ and one of $B$ react (at reaction rate $\kappa$) to form three units of $A$ and one of $C$.  The concentrations $c_{A},$ $c_{B}$, and $c_{C}$ will change in time as the reaction occurs.  Under the assumption of {\em mass-action kinetics}, species $A$ and $B$ react at a rate proportional to the product of their concentrations, where the proportionality constant is the rate constant $\kappa$.  Noting that the reaction yields a net change of two in the amount of $A$, we obtain the first differential equation in the following system:
\begin{align*}
\frac{d}{dt}c_{A}~&=~2\kappa c_{A}c_{B}~, \\
\frac{d}{dt}c_{B} ~&=~-\kappa c_{A}c_{B}~, \\
\frac{d}{dt} c_{C}~&=~\kappa c_{A}c_{B}~.
\end{align*}
The other two equations arise similarly.  Next we include the reverse reaction and switch from additive to multiplicative notation to highlight the monomials that appear in our differential equations; the {\em chemical reaction networks} in this paper will appear with the following notation:
\[
 \begin{xy}<15mm,0cm>:
 (1,0) ="3A+C" *+!L{c_{A}^{3} c_{C} } *{};
 (0,0) ="A+B" *+!R{c_{A}c_{B}} *{};
 (0.55,.12)="k" *+!D{\kappa} *{};
 (0.55, -.5)="k21" *+!D{\kappa '} *{};
   {\ar "A+B"+(0,0.05)*{};"3A+C"+(0,0.05)*{}};
   {\ar "3A+C"+(0,-0.15)*{};"A+B"+(0,-0.15)*{}};
    \end{xy}\]

This network defines differential equations that are each a sum of the monomial contribution from the reactant of each chemical reaction in the network:
\begin{align*}
\frac{d}{dt}c_{A} ~&=~2\kappa c_{A}c_{B}- 2 {\kappa'} c_{A}^{3}c_{C}~, \\
\frac{d}{dt} c_{B}~&=~-\kappa c_{A}c_{B}+ {\kappa'} c_{A}^{3}c_{C}~,\\
\frac{d}{dt} c_{C}~&=~\kappa c_{A}c_{B}- {\kappa'} c_{A}^{3}c_{C}~.
\end{align*}
The recipe for obtaining these differential equations from a chemical reaction network easily generalizes from this example.  However, in order to display the linearity hidden in these non-linear equations, the equations will appear in a different but equivalent form in (\ref{CRN}) below.

We now establish the notation for this paper, following \cite{TDS}.  A {\em chemical reaction network} is a finite directed graph 
whose vertices are labeled by monomials and whose edges are labeled by parameters. 
Specifically, the digraph is denoted $G = (V,E)$, with vertex set $V = \{1,2,\ldots,n\}$
and edge set $\,E \subseteq \{(i,j) \in V \times V : \,i\not= j \}$.
The vertex $i$ of $G$ represents the $i$th chemical complex and is labeled by the monomial
$$ c^{y_i} \,\,\, = \,\,\, c_1^{y_{i1}} c_2^{y_{i2}} \cdots  c_s^{y_{is}}~. $$
This yields $Y=(y_{ij})$, an $n \times s$-matrix of non-negative integers.
The unknowns $c_1,c_2,\ldots,c_s$ represent the
concentrations of the $s$ species in the network,
and we regard them as functions $c_i(t)$ of time $t$.
The monomial labels form the entries in the following row vector: 
$$ \Psi(c) \quad = \quad \bigl( c^{y_1}, ~ c^{y_2} , ~ \ldots,  ~ c^{y_n} \bigr)~. $$
A network is said to be {\em mass-preserving} if all monomials $c^{y_i}$ have the same degree.  
Each directed edge $(i,j) \in E$ is labeled by a positive parameter $\kappa_{ij}$ which
represents the rate constant in the reaction from 
the $i$-th chemical complex to the $j$-th chemical complex.
A network is {\em reversible} if the graph $G$ is undirected, in which case each undirected edge has two labels $\kappa_{ij}$ and $\kappa_{ji}$.
Let $A_\kappa$ denote the negative of the {\em Laplacian} of
the digraph $G$. In other words $A_\kappa$ is the $n \times n$-matrix
whose off-diagonal entries
are the $\kappa_{ij}$ and whose row sums are zero.
Mass-action kinetics specified by the digraph $G$ is
the dynamical system defined by
\begin{equation}
\label{CRN}
 \frac{d c}{dt} \quad = \quad \Psi(c) \cdot A_\kappa \cdot Y ~. 
 \end{equation}
By decomposing the mass-action equations in this way, we see that they are linear in the $\kappa_{ij}$ by way of the matrix $A_{\kappa}$.
A {\em steady state} (or equilibrium) is a positive concentration vector $c\in \mathbb{R}^s_{>0}$ at which the equations (\ref{CRN}) vanish.  These equations remain in the (stoichiometric) subspace $S$ spanned by the vectors $y_i-y_j$ (where $(i,j)$ is an edge of $G$).  In the earlier example, $y_i-y_j =(-2,-1,1)$, meaning that whenever a reaction occurs, two units of $A$ and one of $B$ are lost, while one unit of $C$ is formed (or vice-versa). Therefore, a trajectory $c(t)$ beginning at a positive vector $c(0)=c^0$ remains in the {\em invariant polyhedron} $P:=(c^0+S) \cap \mathbb{R}^s_{\geq 0}$.  {\em Multistationarity} refers to the existence of more than one steady state in some invariant polyhedron. A chemical reaction network may admit multistationarity for all, some, or no choices of positive parameters $\kappa_{ij}$.

Horn initiated the investigation of small chemical reaction networks by enumerating networks 
comprised of ``short complexes,'' those whose corresponding monomials $c^{y} $ have degree at most two \cite{HornCGph,Horn73}. 
Networks that consist of at most three short complexes do not permit multiple steady states.

The next section establishes that the following graph, which we call the {\em Square}, depicts a smallest reversible multistationary chemical reaction network:
\[
 \begin{xy}<15mm,0cm>:
 (0,1) ="c1^3" *+!DR{c_1^3} *{};
 (1,1) ="c1c2^2" *+!DL{c_1 c_2^2} *{};
 (1,0) ="c2^3" *+!UL{c_2^3} *{};
 (0,0) ="c1^2c2" *+!UR{c_1^2 c_2} *{};
 (0.5,1.18)="k12" *+!D{\kappa_{12}} *{};
 (0.5, 0.8)="k21" *+!D{\kappa_{21}} *{};
 (1.21,.5)="k23" *+!L{\kappa_{23}}*{};
 (0.6,.5)="k32" *+!L{\kappa_{32}}*{};
 (.5,-.3)="k34" *+U{\kappa_{34}}*{};
 (.5,.1)="k43" *+U{\kappa_{43}}*{};
 (-.33,.5)="k41" *+R{\kappa_{41}}*{};
 (.2,.5)="k14" *+R{\kappa_{14}}*{};
   {\ar "c1^2c2"+(-0.15,0)*{};"c1^3"+(-0.15,0)*{}};
   {\ar "c2^3"+(0.05,0)*{};"c1c2^2"+(0.05,0)*{}};
   {\ar "c1^3"+(-0.05,0)*{};"c1^2c2"+(-0.05,0)*{}};
   {\ar "c1c2^2"+(+0.15,0)*{};"c2^3"+(+0.15,0)*{}};
   {\ar "c1^3"+(0,+0.15)*{};"c1c2^2"+(0,+0.15)*{}};
   {\ar "c1^2c2"+(0,-0.05)*{};"c2^3"+(0,-0.05)*{}};
   {\ar "c1c2^2"+(0,0.05)*{};"c1^3"+(0,0.05)*{}};
   {\ar "c2^3"+(0,-0.15)*{};"c1^2c2"+(0,-0.15)*{}};
    \end{xy}\]
In the horizontal reactions, two units of species one are transformed into two of species two (or vice-versa), while a third unit remains unchanged by the reaction.  In the vertical reactions, only one is transformed.

The Square appeared in non-reversible form as networks 7-3 in \cite{GenMAK} and 4.2 in \cite{FeinOsc}.  The matrices whose product defines the dynamical system (\ref{CRN}) follow:
\[
\Psi(c) \quad 		= \quad \bigl( c_1^3, ~ c_1c_2^2 , ~ c_2^3, ~ c_1^2c_2 \bigr)~, \]
\[
A_{\kappa} \quad	= \quad  \begin{pmatrix}
 - \kappa_{12} - \kappa_{14} & \kappa_{12} & 0 & \kappa_{14} \\
  \kappa_{21} & -\kappa_{21} - \kappa_{23} & \kappa_{23} & 0 \\
  0 & \kappa_{32} & -\kappa_{32} - \kappa_{34} & \kappa_{34} \\
  \kappa_{41} & 0 & \kappa_{43} & -\kappa_{41} - \kappa_{43} 
  \end{pmatrix}~, \]
\[
Y \quad 			= \quad \begin{pmatrix}
3 & 0 \\
1 & 2 \\
0 & 3 \\
2 & 1 \\
\end{pmatrix} ~.  
\]
There may be two or even three steady states in each invariant polyhedron $P$; Example \ref{ex1} in the next section provides a choice of positive rate constants $\kappa_{ij}$ that give rise to three steady states.  Sections 4 and 5 determine precisely which parameters give rise to two steady states and which yield three.  Moreover, we compute this semi-algebraic parametrization for all networks on the same four vertices as the Square, in other words, networks with complexes $c_1^3, ~ c_1c_2^2 , ~ c_2^3,$ and $c_1^2c_2$.  The parametrization is captured in Table 1 and can be computed ``by hand,'' but larger systems may require  techniques of computational real algebraic geometry \cite{BPR}.  For example, our problem of classifying parameters according to number of steady states is labeled as Problem P2  in \cite{StabA}, where it
is addressed with computer algebra methods.

\section{The Smallest Multistationary Network} 

Following equation (7) of \cite{TDS}, the Matrix-Tree Theorem defines the following polynomials in the rate constants of the Square:
\begin{align*}
K_1 ~=~ \kappa_{23} \kappa_{34} \kappa_{41} +\kappa_{21} \kappa_{34} \kappa_{41}+\kappa_{21} \kappa_{32}\kappa_{41} +\kappa_{21} \kappa_{32} \kappa_{43}~, \\
K_2 ~=~ \kappa_{14 } \kappa_{32 } \kappa_{43 } +\kappa_{12} \kappa_{34} \kappa_{41}+\kappa_{12} \kappa_{32} \kappa_{41}+\kappa_{12} \kappa_{32} \kappa_{43}~, \\
K_ 3 ~=  ~ \kappa_{ 14} \kappa_{23 } \kappa_{ 43} +\kappa_{14} \kappa_{21} \kappa_{43}+\kappa_{12} \kappa_{23} \kappa_{41}+\kappa_{12} \kappa_{23} \kappa_{43}~, \\
K_ 4 ~=~ \kappa_{ 14} \kappa_{23 } \kappa_{34 } +\kappa_{14} \kappa_{21} \kappa_{34}+\kappa_{14} \kappa_{21} \kappa_{32}+\kappa_{12} \kappa_{23} \kappa_{34} ~.
\end{align*}
Theorem 7 of \cite{TDS} provides an ideal $M_G$ that is toric in these $K_i$ coordinates, and the variety of $M_G$ is the moduli space of toric dynamical systems on the Square.  In this case, the ideal $M_G$ is the {\em twisted cubic curve} in the $K_i$ coordinates, generated by the 2$\times$2-minors of the following matrix:
\begin{align}\label{twistedMatrix}
\begin{pmatrix}
K_1 \quad & K_2 \quad &K_4 \\
K_4 \quad &  K_3 \quad & K_2
\end{pmatrix}~.
\end{align}
Theorem 7 of \cite{TDS} says that for a given choice of positive rate constants $\kappa_{ij}$, the equations (\ref{CRN}) define a toric dynamical system if and only if the minors of the matrix (\ref{twistedMatrix}) vanish.  In general the codimension of $M_G$ is the {\em deficiency} of a network; see Theorem 9 of \cite{TDS}.  Here the deficiency is two.  Recall that a {\em toric dynamical system} is a dynamical system (\ref{CRN})
for which the algebraic equations
  $\,\Psi(c) \cdot A_\kappa  = 0\,$
  admit a strictly positive solution
  $c^* \in \mathbb{R}^s_{>0}$; this solution is called a complex balancing steady state \cite{GenMAK}.  In this case there is a unique steady state in each invariant polyhedron $P$, so multistationarity is ruled out.  Toric dynamical systems exhibit further nice properties; for details, see \cite{TDS,FeinLectures,GenMAK}.

It is no coincidence that the original monomials of the Square, namely $ c_1^3,$ $c_1c_2^2,$ $c_2^3,$ $c_1^2c_2$, parametrize the twisted cubic curve.  In fact, the following general result follows from Theorem 9 in \cite{TDS}.
\begin{proposition}
Assume that a chemical reaction network $G$ is strongly connected and all monomials $c^{y_{i}}$ have the same total degree.  Then the toric variety parametrized by $\Psi(c)$ coincides with the variety of $M_{G}$. 
\end{proposition}
For the Square, each one-dimensional invariant polyhedron $P$ is defined by some positive concentration total $T=c_1+c_2$.   The steady states in $P$ correspond precisely to the positive roots of the following cubic polynomial:
$$p_S(x)~=~(-2 \kappa_{12}-\kappa_{14})x^3+ (\kappa_{41}-\kappa_{43})x^2+(2\kappa_{21}-\kappa_{23})x+(\kappa_{32}+2\kappa_{34}) ~ ;
$$
this polynomial arises by substituting $x:=c_1/{c_2}$ into the equation $d{c_1}/{dt}=-d{c_2}/{dt}$.
From this point of view, we reach some immediate conclusions.  First, the {\em algebraic degree} of this system is three, which bounds the number of steady states.  Second, the number of steady states and their stability depend only on the rate parameters $\kappa_{ij}$, and not on the invariant polyhedron $P$ or equivalently the choice of total concentration $T$.
Also, by noting that $p_S(x)$ is positive at $x=0$ and is negative for large $x$, we see that the Square admits at least one steady state for any choice of rate constants. 
Recall that the discriminant of a univariate polynomial $f$ is a polynomial that vanishes precisely when $f$ has a multiple root over the complex numbers \cite{PolyBk}.   {\tt Maple} computes the discriminant of $p_S$ to be the following polynomial:
\begin{align*}
&   - 108 \kappa_{12}^2 \kappa_{32}^2 
 - 432 \kappa_{12}^2 \kappa_{32} \kappa_{34} 
 - 432 \kappa_{12}^2 \kappa_{34}^2 
 - 108 \kappa_{12} \kappa_{14} \kappa_{32}^2 
\\ & - 432 \kappa_{12} \kappa_{14} \kappa_{32} \kappa_{34}
  -432 \kappa_{12} \kappa_{14} \kappa_{34}^2 
  + 64 \kappa_{12} \kappa_{21}^3 
- 96 \kappa_{12} \kappa_{21}^2 \kappa_{23} 
+ 48 \kappa_{12} \kappa_{21} \kappa_{23}^2 
\\ & - 72 \kappa_{12} \kappa_{21} \kappa_{32} \kappa_{41}
 + 144 \kappa_{12} \kappa_{21} \kappa_{32} \kappa_{43}
  - 144 \kappa_{12} \kappa_{21} \kappa_{34} \kappa_{41} 
+ 288 \kappa_{12} \kappa_{21} \kappa_{34} \kappa_{43} 
\\ & - 8 \kappa_{12} \kappa_{23}^3 
+ 36 \kappa_{12} \kappa_{23} \kappa_{32} \kappa_{41} 
- 72 \kappa_{12} \kappa_{23} \kappa_{32} \kappa_{43} 
+ 72 \kappa_{12} \kappa_{23} \kappa_{34} \kappa_{41} 
\\ & - 144 \kappa_{12} \kappa_{23} \kappa_{34} \kappa_{43}
 - 27 \kappa_{14}^2 \kappa_{32}^2 
- 108 \kappa_{14}^2 \kappa_{32} \kappa_{34} 
 - 108 \kappa_{14}^2 \kappa_{34}^2 
 + 32 \kappa_{14} \kappa_{21}^3 
\\ & - 48 \kappa_{14} \kappa_{21}^2 \kappa_{23} 
+ 24 \kappa_{14} \kappa_{21} \kappa_{23}^2 
- 36 \kappa_{14} \kappa_{21} \kappa_{32} \kappa_{41} 
 + 72 \kappa_{14} \kappa_{21} \kappa_{32} \kappa_{43}
\\ &  - 72 \kappa_{14} \kappa_{21} \kappa_{34} \kappa_{41} 
+ 144 \kappa_{14} \kappa_{21} \kappa_{34} \kappa_{43} 
- 4 \kappa_{14} \kappa_{23}^3 
 + 18 \kappa_{14} \kappa_{23} \kappa_{32} \kappa_{41} 
\\ &  - 36 \kappa_{14} \kappa_{23} \kappa_{32} \kappa_{43} 
+ 36 \kappa_{14} \kappa_{23} \kappa_{34} \kappa_{41} 
- 72 \kappa_{14} \kappa_{23} \kappa_{34} \kappa_{43}
  + 4 \kappa_{21}^2 \kappa_{41}^2 
\\ &  - 16 \kappa_{21}^2 \kappa_{41} \kappa_{43} 
+ 16 \kappa_{21}^2 \kappa_{43}^2 
- 4 \kappa_{21} \kappa_{23} \kappa_{41}^2
 + 16 \kappa_{21} \kappa_{23} \kappa_{41} \kappa_{43}
  - 16 \kappa_{21} \kappa_{23} \kappa_{43}^2 
\\ & + \kappa_{23}^2 \kappa_{41}^2 
 - 4 \kappa_{23}^2 \kappa_{41} \kappa_{43} 
+ 4 \kappa_{23}^2 \kappa_{43}^2 
 - 4 \kappa_{32} \kappa_{41}^3 
+ 24 \kappa_{32} \kappa_{41}^2 \kappa_{43} 
 - 48 \kappa_{32} \kappa_{41} \kappa_{43}^2 
\\ & + 32 \kappa_{32} \kappa_{43}^3 
 - 8 \kappa_{34} \kappa_{41}^3 
+ 48 \kappa_{34} \kappa_{41}^2 \kappa_{43}
 - 96 \kappa_{34} \kappa_{41} \kappa_{43}^2 
 + 64 \kappa_{34} \kappa_{43}^3 ~.
\end{align*}
As $p_S$ is cubic and has at least one positive root, its discriminant is negative if and only if $p_S$ has one real root and one pair of complex conjugate roots; in this case, the Square has precisely one steady state.  When the discriminant is non-negative, the system may admit one, two, or three steady states; we analyze this case fully in the next section.  
\begin{example}\label{ex1}
Consider the following rate constants for the Square:
$$ (\kappa_{12}, ~\kappa_{14}, ~\kappa_{21}, ~\kappa_{23}, ~\kappa_{32}, ~\kappa_{34}, ~\kappa_{41}, ~\kappa_{43} ) \quad = \quad (1/4 ,~ 1/2 ,~1   ,~  13 ,~ 1 ,~   2,~ 8,~1) ~.$$
This yields $p_S (x)=-x^3+6x^2-11+6$, which has three positive roots: $x=1$, 2, and 3.  This is an instance of bistability; it is easy to determine that $x=1$ and $x=3$ correspond to {\em stable} steady states, while the third is unstable.  In the next section we determine the conditions for an arbitrary vector of rate constants to admit one, two, or three steady states.
\end{example}

Recalling the definitions given earlier, the Square has the following properties: 
the number of complexes is $n=4$,
the number of connected components of $G$ is $l=1$,
the number of species is $s=2$, and
the dimension of any invariant polyhedron is $\sigma= 1$. The main result of this section states that this network is minimal with respect to each of these four parameters.

\begin{theorem}\label{square}
The Square is a smallest multistationary, mass-preserving,
reversible chemical reaction network with respect to each of the following parameters:
the number of complexes,
the number of connected components of $G$,
the number of species, and
the dimension of an invariant polyhedron.
\end{theorem}

\begin{proof}
First $l=1$ and $\sigma=1$ are clearly minimal.  Next any mass-preserving system with $n\leq 2$ or $s=1$ has no reactions or has deficiency zero.  Finally, an $n=3$ system has deficiency zero or one; in the
deficiency one case, the Deficiency One Theorem of Feinberg rules out the possibility of
multistationarity \cite{FeinDefOne}.
\qed
\end{proof}

Among all mass-preserving multistationary systems that share these four minimal parameters, the Square is distinguished because its monomials are of minimal degree.  A connected network of lower degree would consist of at most three of Horn's ``short'' complexes \cite{HornCGph}.  

We now discuss the possible connection of the Square to biology by comparing it to the following simple 
network:
\begin{align}
\label{memory}
c_x c_y  & \leftrightarrows c_y^2  \\
c_x & \leftrightarrows c_y . \notag
\end{align}
Network (\ref{memory}) is a modified version of the following molecular switch mechanism proposed by Lisman \cite{Lisman}: 
\begin{eqnarray*}
c_x c_y  & \leftrightarrows c_{xy} & \longrightarrow c_y^2  \\
c_y c_p & \leftrightarrows c_{yp} & \longrightarrow  c_x c_p .
\end{eqnarray*}
Here $x$ denotes a kinase in an inactive state, $y$ is the active version, and $p$ is a phosphatase.  In the first reactions, $y$ catalyzes the phosphorylation of $x$, turning $x$ into $y$; the second reactions correspond to dephosphorylation.  By skipping the binding steps, making all reactions reversible, and noting that removing $p$ effectively scales the second reaction rate constant, we obtain the network (\ref{memory}).  The reactions of (\ref{memory}) are similar to $c_1^2 c_y   \leftrightarrows c_2^3$ and $c_1^3  \leftrightarrows c_2^3$, which are reactions in the generalization of the Square network examined in the next section; this suggests the possible biological relevance of the reactions of the Square.  For example $c_1^2 c_2 \longrightarrow c_2^3$ can be viewed as a reaction in which species two catalyzes the reaction $c_1^2 \longrightarrow c_2^2$.  Such a positive feedback loop--in which a high amount of some species $y$ encourages the further production of the same species--occurs in biological settings.  For example, the recent work of Dentin {\em et al.} finds that high glucose levels in diabetic mice promote further glucose production in the liver, which is triggered by the binding of glucose (which we may view as $y$) to the transcription factor CREB ($x$) \cite{Dentin}.

This paper focuses on the Square and more generally, the networks that share the same complexes as the Square.  In the following section, we shall determine which of these are bistable.  
The one with the fewest edges is the only one with two connected components rather than one, and is featured in the last section.

\section{Parametrizing Multistationarity}

The aim of this section is similar to that of Conradi {\em et al.} \cite{MAPK}, which determined the full set of parameters that give rise to multistationarity for a biochemical model describing a single layer of a MAPK cascade.  
However we additionally determine the precise number of steady states: zero, one, two, or three, and determine their stability.  The family of networks we consider are those that have the same four complexes as the Square.  In other words, we classify subnetworks of the complete network depicted here:
\[
 \begin{xy}<15mm,0cm>:
 (0,1) ="c1^3" *+!DR{c_1^3} *{};
 (1,1) ="c1c2^2" *+!DL{c_1 c_2^2} *{};
 (1,0) ="c2^3" *+!UL{c_2^3} *{};
 (0,0) ="c1^2c2" *+!UR{c_1^2 c_2} *{};
   {\ar "c1^2c2"+(-0.15,0)*{};"c1^3"+(-0.15,0)*{}};
   {\ar "c2^3"+(0.05,0)*{};"c1c2^2"+(0.05,0)*{}};
   {\ar "c1^3"+(-0.05,0)*{};"c1^2c2"+(-0.05,0)*{}};
   {\ar "c1c2^2"+(+0.15,0)*{};"c2^3"+(+0.15,0)*{}};
   {\ar "c1^3"+(0,+0.15)*{};"c1c2^2"+(0,+0.15)*{}};
   {\ar "c1^2c2"+(0,-0.05)*{};"c2^3"+(0,-0.05)*{}};
   {\ar "c1c2^2"+(0,0.05)*{};"c1^3"+(0,0.05)*{}};
   {\ar "c2^3"+(0,-0.15)*{};"c1^2c2"+(0,-0.15)*{}};
   {\ar "c1^3"+(0.05,-0.05)*{};"c2^3"+(0.0,0.05)*{}};
   {\ar "c2^3"+(-0.05,0)*{};"c1^3"+(0.0,-.1)*{}};
   {\ar "c1^2c2"+(0.05,0.05)*{};"c1c2^2"+(-.05,-.05)*{}};
   {\ar "c1c2^2"+(0,-.10)*{};"c1^2c2"+(0.1,0.00)*{}};
    \end{xy}\]
\\
Each of the twelve rate constants $\kappa_{ij}$ is permitted to be zero, 
which defines the parameter space $\mathbb{R}_{\geq 0}^{12}$ of  dynamical systems.
The main result of this section is summarized in Table 1, which is the semi-algebraic decomposition of the twelve-dimensional parameter space according to the number of steady states of the dynamical system.  The conditions listed there make use of certain polynomials in the rate constants, including the signed coefficients of the polynomial $p$: 
\begin{align*}
S_0 \quad &=\quad    2\kappa_{12} + 3\kappa_{13}+\kappa_{14}  ~, \\ 
S_1 \quad &= \quad  \kappa_{41}-\kappa_{42}-2\kappa_{43} ~, \\ 
S_2 \quad &= \quad  -2\kappa_{21}+\kappa_{23} -\kappa_{24}~,  \\ 
S_3 \quad &= \quad   3\kappa_{31}+\kappa_{32}+2\kappa_{34} ~,  
\end{align*}
where $p$ generalizes the polynomial $p_S$ from the Square:
\begin{align*}\label{p}
p(x)\quad =\quad	-S_0x^3 
				+S_1 x^2 
				 -S_2x
				+S_3 ~ .
\end{align*} 
\begin{table} \label{table1}
\begin{center} 
\caption{Classification of dynamical systems arising from non-trivial (having at least one reaction) networks with complexes $c_1^3, ~ c_1c_2^2 , ~ c_2^3, ~ c_1^2c_2$.  Listed are the number of steady states and the number of steady states that are stable.  The discriminant of $p$ is denoted by $D$.  The signed coefficients of $p$ are denoted by $S_0$, $S_1$, $S_2$, and $S_3$.  The triple root condition consists of the equations (4). }
 \begin{tabular}{ | c || c | c |}
   \hline
    Condition  & Steady states & Stable states   \\ \hline \hline

	$D<0$ and $S_0 S_3=0$						& 0 & 0 \\
	$D<0$ and {\em else}						& 1 & 1 \\  \hline

	$D>0$ and $S_0,S_1,S_2,S_3>0$			& 3 & 2 \\
	$D>0$ and $S_0,S_1,S_2>0$ and $S_3=0$				& 2 & 1 \\
	$D>0$ and $S_1,S_2,S_3>0$ and $S_0=0$				& 2 & 1 \\
	$D>0$ and $S_0=S_3=0$ and $S_1 S_2<0$	& 0 & 0 \\
	$D>0$ and {\em else}						& 1 & 1 \\ \hline
	
	$D=0$ and $S_0,S_1,S_2,S_3>0$ and triple root condition				& 1 & 1 \\
	$D=0$ and $S_0,S_1,S_2,S_3>0$ without triple root condition			& 2 & 1 \\
	$D=0$ and $S_1\leq S_0=0 \leq S_2$	and $S_3>0$		& 0 & 0 \\
	$D=0$ and $S_1\leq S_3=0 \leq S_2$ and $S_0>0$			& 0 & 0 \\
	$D=0$ and {\em else}					& 2 & 1 \\
\hline
   \hline
 \end{tabular}
\end{center}
\end{table}

We now derive the entries of Table 1 for those networks without boundary steady states (this includes the case of the Square).  These cases are precisely the ones in which $S_0>0$ and $S_3>0$.  Our approach is simply to determine the conditions on the coefficients of $p$ for the polynomial to have one, two, or three positive roots. 

In this twelve-parameter case, the discriminant of $p$ is a homogeneous degree-four polynomial with 113 terms.  For the same reason as that for the Square, there is one steady state when the discriminant is negative. Now assume that the discriminant is non-negative.   Then $p$ has three real roots, counting multiplicity; recall that the positive ones correspond to the steady states of the chemical reaction network.  Now the constant term of a monic cubic polynomial is the negative of the product of its roots, so by examining the sign of the constant term of $p$, we conclude that $p$ has either one positive root and two negative roots, or three positive roots.  Continuing the sign analysis with the other coefficients of $p$, we conclude that there are three positive roots if and only if $S_1>0$ and $S_2>0$.  We proceed by distinguishing between the cases when the discriminant is positive or zero.  If the discriminant is positive, then we have derived criteria for having one or three steady states; this is because the roots of $p$ are distinct.  If the discriminant is zero, then in the case of one positive root, the two negative roots come together (one steady state).  In the case of discriminant zero and three positive roots, then at least two roots come together (at most two steady states); a triple root occurs if and only if the following {\em triple root condition} holds:
\begin{align} 
3 S_0 S_2 ~= ~S_1^2 \quad \text{and} \quad 27 S_0^2S_3~=~S_1^3~.
\end{align}
These equations are precisely what must hold in order for $p$ to have the form $p(x)=-(x-\alpha)^3$.  Finally, stability analysis in this one-dimensional system is easy, and this completes the analysis for the networks without boundary steady states.  The remaining cases can be classified similarly to complete the entries of Table 1. 
To parametrize the behavior of the Square, we simply reduce to the case when each of its parameters $\kappa_{12},~\kappa_{14}, ~\kappa_{21}, ~\kappa_{23}, ~\kappa_{32}, ~\kappa_{34}, ~\kappa_{41},$ and $\kappa_{43} $ are positive and all others are zero.  

By determining which sign vectors in $(0,+)^{12}$ can be realized by a vector of parameters that yields multistationarity, we find a necessary and sufficient condition for a directed graph on the four complexes of the Square to admit multistationarity.  This condition is that the graph must include the edges labeled by $\kappa_{23}$ and $\kappa_{41}$ and at least one edge directed from the vertex $c_1^3$ or $c_2^3$.  
In this case, for appropriate rate parameters arising from Table 1, the dynamical system has multiple steady states.
 Therefore, we can enumerate the reversible networks on the four complexes that admit multistationarity: there is one network with all six (bi-directional) edges, four with five edges, six (including the Square) with four edges, four with three edges, and one with two edges.  These sixteen networks comprise the family of ``smallest'' multistationary networks.  For the two-edge network, the decomposition from Table 1 is depicted in Figure~\ref{VertFig} in the next section.
 
\section{Subnetworks of the Square}

Subnetworks of the Square are obtained by removing edges.  From the parametrization in the previous section, we know that up to symmetry between $c_1$ and $c_2$, only two reversible subnetworks of the Square exhibit multiple steady states. 

The first network is obtained by removing the bottom edge of the Square. 
In other words $A_{\kappa}$ is replaced by
\[
A_{\kappa} \quad	= 	\quad  \begin{pmatrix}
 - \kappa_{12} - \kappa_{14} & \kappa_{12} & 0 & \kappa_{14} \\
  \kappa_{21} & -\kappa_{21} - \kappa_{23} & \kappa_{23} & 0 \\
  0 & \kappa_{32} & -\kappa_{32}  & 0 \\
  \kappa_{41} & 0 & 0 & -\kappa_{41} 
  \end{pmatrix}~. \]

In this subnetwork, the four parameters of Theorem \ref{square} are the same as those of the Square.  The system is a toric dynamical system if and only if the following four binomial generators of $M_G$ vanish:
\begin{align*}
\kappa_{14} \kappa_{32}  ~&-~\kappa_{23}  \kappa_{41} ~ , \\ 
 \kappa_{12}  \kappa_{32}  \kappa_{41} ~&-~\kappa_{14}  \kappa_{21}  \kappa_{23}  ~, \\ 
\kappa_{14}^2 \kappa_{21}   ~&-~ \kappa_{12}  \kappa_{41}^2~,\\ 
 \kappa_{12} \kappa_{32}^2 ~&-~ \kappa_{21} \kappa_{23}^2 ~. 
\end{align*}
We note that both $\kappa_{23}$ times the third binomial and $\kappa_{14}$ times the fourth binomial are in the ideal generated by the first two binomials.  Therefore, an assignment of positive parameters for this network defines a toric  dynamical system if and only if the following two equations hold: $\kappa_{14} \kappa_{32}  ~=~\kappa_{23}  \kappa_{41}$ and
$ \kappa_{12}  \kappa_{32}  \kappa_{41} ~=~ \kappa_{14}  \kappa_{21}  \kappa_{23}
$.

\begin{figure}[h]
\centering
\includegraphics[scale=0.7]{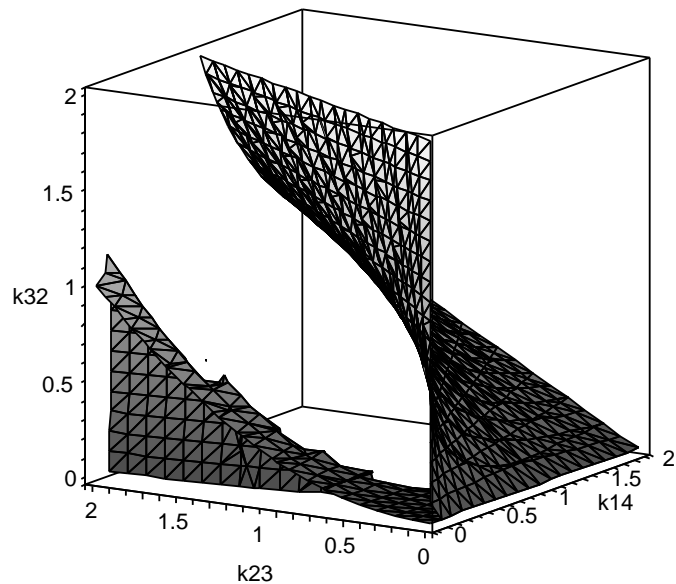}
\caption{This depicts the semi-algebraic decomposition of Section 4 for the subnetwork of the Square in which only the vertical edges remain and $\kappa_{41}=1$.  At the left is the discriminant-zero locus.  Parameter vectors lying below this surface give rise to dynamical systems with three steady states.  Those above the surface yield one steady state; these include parameters of the toric dynamical systems, which are the points on the Segre variety which appears on the right.  Parameters on the discriminant-zero locus correspond to systems with either one (if $3 \kappa_{14} \kappa_{32} = \kappa_{23}$) or two steady states.  This figure was created using {\tt Maple}. }
\label{VertFig}
\end{figure}

The second subnetwork of the Square is obtained by removing one additional edge, the one between the vertices labeled by $c_1^3$ and $c_1c_2^2$.  The new $A_{\kappa}$ is
\[
A_{\kappa} \quad	=	\quad  \begin{pmatrix}
 - \kappa_{14} & 0 & 0 & \kappa_{14} \\
  0&  - \kappa_{23} & \kappa_{23} & 0 \\
  0 & \kappa_{32} & -\kappa_{32}  & 0 \\
  \kappa_{41} & 0 & 0 & -\kappa_{41} 
  \end{pmatrix}~. \]
The network graph $G$ is now disconnected, and $p$ reduces to

$$
p(x) ~=~ -\kappa_{14}  x^3+\kappa_{41}  x^2-\kappa_{23}  x+\kappa_{32} ~ .
$$
The discriminant of $p$ is 
$$
D \quad = \quad -27 \kappa_{14}^2 \kappa_{32}^2-4 \kappa_{14} \kappa_{23}^3+18 \kappa_{14} \kappa_{23} \kappa_{32} \kappa_{41}+\kappa_{23}^2 \kappa_{41}^2-4 \kappa_{32} \kappa_{41}^3
 ~.
$$
Further, the toric condition reduces to the single equation
$$
\kappa_{23}\kappa_{41} ~=~ \kappa_{14} \kappa_{32}~ ,$$
which defines the {\em Segre variety}.  A single equation suffices to define the space of toric dynamical systems; this corresponds to the fact that this subnetwork has deficiency one, while the previous subnetwork has deficiency two.  The semi-algebraic decomposition of the previous section for this four-parameter network can be depicted in three dimensions by setting one parameter to be one, in other words, by scaling the equations (\ref{CRN}); this is displayed in Figure~\ref{VertFig}.

We remark that Horn and Jackson performed the same parametrization for the following special rate constants:
$$ (\kappa_{12}, ~\kappa_{14}, ~\kappa_{21}, ~\kappa_{23}, ~\kappa_{32}, ~\kappa_{34}, ~\kappa_{41}, ~\kappa_{43} ) \quad = \quad (\epsilon  ,~ 0 ,~1   ,~  0 ,~ \epsilon ,~  0,~ 1,~0) ~, $$
where $\epsilon >0$.  Their results are summarized as Table 1 in {\cite{GenMAK}.  Their analysis notes that any instance of three steady states can be lifted to establish the same in the (reversible) Square.  In other words, in a small neighborhood in $\mathbb{R}_{\geq 0}^8$ of a vector of parameters that yields three steady states of the directed Square, there is a vector of parameters for the bi-directional Square that also exhibits multistationarity.  
The specific criterion for when lifting of this form is possible appears in Theorem 2 of 
Conradi {\em et al.} \cite{subnetwork}.  As this approach is widely applicable, further analysis of small networks may be fruitful for illuminating the dynamics of larger biochemical networks.  

We have seen that the family of Square networks is the smallest class of bistable mass-action kinetics networks.  Whether nature has implemented one of these (perhaps with additional components to provide robustness) in a biological setting is as yet unknown, but it is also remarkable that these networks exhibit a simple switch mechanism, which we now explain.  Consider the case of three steady states.  The corresponding positive roots $x_{1}<x_{2}<x_{3}$ of $p$ in Section 4 are the equilibria for the ratio of concentrations $c_{1}/c_{2}$.  To switch from the low stable equilibrium $x_{1}$ to the high stable equilibrium $x_{3}$ is easy: simply increase the concentration ratio $c_{1}/c_{2}$ past $x_{2}$, and the dynamics will do the rest.

\subsubsection*{Acknowledgments}  Bernd Sturmfels posed the question of determining the smallest bistable network and provided guidance for this work.
We thank Carsten Conradi, Gheorge Craciun, Lior Pachter, and J\"{o}rg Stelling for helpful discussions.
Anne Shiu was supported by a Lucent Technologies Bell Labs Graduate Research Fellowship.  

%
%

%
\end{document}